\documentclass[12pt]{article}
\usepackage{ifpdf}
\usepackage{graphicx,subfigure}
\usepackage[english]{babel}
\usepackage[latin1]{inputenc}
\usepackage{mathrsfs}
\usepackage{amsfonts}
\usepackage{amsmath}

\newcommand{\diff}{\mathrm{d}}
\newtheorem{mydef}{Definition}[section]

\ifx\pdfoutput\undefined
\usepackage[hypertex]{hyperref}
\else
\usepackage[pdftex,hypertexnames=false]{hyperref}
\fi

\begin{document}
\thispagestyle{empty}
\vfill
\begin{center} {\Large Model-free control of non-minimum phase systems and switched systems} \end{center}
\vfill
\begin{center}{\sc Lo\"ic MICHEL} \end{center}
\begin{center}{ \it {GRÉI - D\'epartement de g\'enie \'electrique et g\'enie informatique \\
Université du Qu\'ebec \`a Trois-Rivi\`eres \\
C.P. 500, Trois-Rivi\`{e}res, G9A 5H7, Canada, QC }} \end{center}
\vfill
\begin{quotation}                
\begin{center} {\bf Abstract} \end{center}
\noindent
This brief\footnote{This work is distributed under CC license \url{http://creativecommons.org/licenses/by-nc-sa/3.0/}} presents a simple
derivation of the standard model-free control for the non-minimum phase systems. The robustness of the proposed method is studied in simulation considering the case of switched systems.
\end{quotation}
\vfill

\vfill
\def\thefootnote{\fnsymbol{footnote}}
\setcounter{footnote}{0}

\newpage
\setcounter{page}{1}
\section{Introduction}
The model-free control methodology, originally proposed by \cite{esta}, has been widely successfully applied to many mechanical
and electrical processes. The model-free control provides good performances in disturbances rejection and an efficient robustness
to the process internal changes. The control of non-minimum phase systems has been deeply studied and successful methods have been
proposed (e.g. \cite{astrom} \cite{Bai} \cite{isidori} \cite{chen} \cite{Benosman} \cite{Gurumoorthy} \cite{Karagiannis} \cite{Barkana}).
Since the model-free control can not {\it a priori}  stabilize a non-minimum phase system \cite{esta}, we propose a possible derivation of the original model-free control law, dedicated to the control of non-minimum phase systems. The dynamic performances are especially tested in the case of switched systems.

The paper is structured as follows. Section II presents an overview of the model-free control methodology including its advantages in comparison
with classical methodologies. Section III discusses the application of the modified model-free control, called NM-model-free control, for non-minimum phase systems. Some concluding remarks may be found in Section IV.

\section{Model-free control: a brief overview}\label{mfc}
\subsection{General principles}
\subsubsection{The ultra-local model}
We only assume that the plant behavior is well approximated in its
operational range by a system of ordinary differential equations,
which might be highly nonlinear and time-varying.\footnote{See
\cite{esta,Fliess} for further details.} The system, which is SISO,
may be therefore described by the input-output equation
\begin{equation}\label{es}
E (t, y, \dot{y}, \dots, y^{(\iota)}, u, \dot{u}, \dots,
u^{(\kappa)}) = 0
\end{equation}
where
\begin{itemize}
\item $u$ and $y$ are the input and output variables,
\item $E$, which might be unknown, is assumed to be a
sufficiently smooth function of its arguments.
\end{itemize}
Assume that for some integer $n$, $0 < n \leq \iota$,
$\frac{\partial E}{\partial y^{(n)}} \not\equiv 0$. From the
implicit function theorem we may write locally
$$
y^{(n)} = \mathfrak{E} (t, y, \dot{y}, \dots, y^{(n - 1)}, y^{(n +
1)}, \dots, y^{(\iota)}, u, \dot{u}, \dots, u^{(\kappa)})
$$
By setting $\mathfrak{E} = F + \alpha u$ we obtain the {\it
ultra-local} model.

\begin{mydef}\label{mydef-modele_F}
\cite{esta} If $u$ and $y$ are respectively the variables of input and output of a system to be controlled,
then this system admits the ultra-local model defined by:
\begin{equation}\label{mod}
y^{(n)} = F + \alpha u
\end{equation}
where
\begin{itemize}
\item $\alpha \in \mathbb{R}$ is a {\em non-physical} constant parameter,
such that $F$ and $\alpha u$ are of the same magnitude;
\item the numerical value of $F$, which contains the whole ``structural information'',
is determined thanks to the knowledge of $u$, $\alpha$, and of the
estimate of the derivative $y^{(n)}$.
\end{itemize}
\end{mydef}
In all the numerous known examples it was possible to set $n = 1$ or
$2$.

\subsubsection{Numerical value of $\alpha$}\label{beta}Let
us emphasize that one only needs to give an approximate numerical
value to $\alpha$. It would be meaningless to refer to a precise
value of this parameter.

\subsection{Intelligent PI controllers}

\subsubsection{Generalities}

\begin{mydef}
\cite{esta} \label{mydef-iPI} we close the
loop via the {\em intelligent PI controller}, or {\em i-PI}
controller,
\begin{equation}\label{eq:ipi}
u = - \frac{F}{\alpha} + \frac{\dot{y}^{(n)\ast} }{\alpha}  + \mathcal{C}(\varepsilon)
\end{equation}
where
\begin{itemize}
\item $y^\ast$ is the output reference trajectory, which
is determined via the rules of flatness-based control
(\cite{flat1,flat2});
\item $e = y^\ast - y$ is the tracking error;
\item $\mathcal{C}(\varepsilon)$ is of the form $K_P \varepsilon + K_I \int \varepsilon$. $K_P$, $K_I$ are the usual tuning gains.
\end{itemize}

Equation (\ref{eq:ipi}) is called model-free control law or model-free law.
\end{mydef}

The i-PI controller \ref{eq:ipi} is compensating the poorly known
term $F$. Controlling the system therefore boils down to the control
of a precise and elementary pure integrator. The tuning of the gains
$K_P$ and $K_I$ becomes therefore quite straightforward.

\subsubsection{Classic controllers}
See \cite{maroc} for a comparison with classic PI controllers.

%
%
%
%
%

\subsection{A first academic example: a stable monovariable linear system}\label{list}
Introduce as in \cite{esta,Fliess} the stable transfer function
\begin{equation}\label{eq:stlin}
\frac{(s+2)^2}{(s+1)^3}
\end{equation}

\subsubsection{A classic PID controller}\label{broida} We apply the
well known method due to Bro\"{\i}da (see, {\it e.g.}, \cite{dinde}) by
approximating System \ref{eq:stlin} via the following delay system
$$
\frac{Ke^{-\tau s}}{(Ts+1)}
$$
$K = 4$, $T = 2.018$, $\tau = 0.2424$ are obtained thanks to
graphical techniques. The gain of the PID controller are then
deduced \cite{dinde}: $K_P =\frac{100(0.4\tau+T)}{120K\tau}=1.8181$,
$K_I = \frac{1}{1.33K\tau}=0.7754$, $K_D = \frac{0.35T}{K}=0.1766$.

\subsubsection{i-PI.}
We are employing $\dot{y} = F + u$ and the i-PI controller
$$
u=-[F]_e+\dot y^\star + \mathcal{C}(\varepsilon)$$ where
\begin{itemize}
\item $[F]_e=[\dot y]_e-u$,
\item $y^\star$ is a reference trajectory,
\item $\varepsilon = y^\star - y$,
\item $\mathcal{C}(\varepsilon)$ is an usual PI controller.
\end{itemize}

\subsubsection{Numerical simulations}

Figure 1(a) shows that the i-PI controller behaves only slightly
better than the classic PID controller (Fig. 1(b)). When taking into account on
the other hand the ageing process and some fault accommodation there
is a dramatic change of situation:
Figure 1(c) indicates a clear cut superiority of our i-PI
controller if the ageing process corresponds to a shift of the pole
from $1$ to $1.5$, and if the previous graphical identification is
not repeated (Fig. 1(d)).

\subsubsection{Some consequences}
\begin{itemize}
\item It might be useless to introduce delay systems of the type
\begin{equation}\label{retard}
T(s) e^{-Ls}, \quad ~ T \in {\mathbb{R}}(s), ~ L \geq 0
\end{equation}
for tuning classic PID controllers, as often done today in spite of
the quite involved identification procedure.
\item This example demonstrates also that the usual mathematical
criteria for robust control become to a large irrelevant.
\item As also shown by this example some fault accommodation may also
be achieved without having recourse to a general theory of
diagnosis.
\end{itemize}

\begin{figure}[!h]
\centering
\subfigure[\footnotesize i-PI control]{\includegraphics[width=2in]{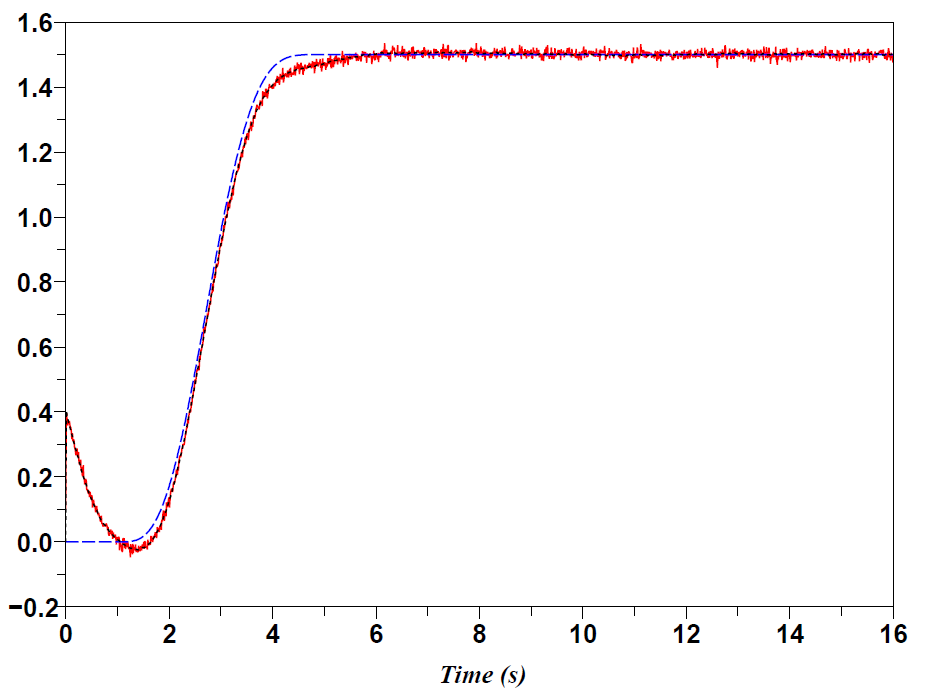}\label{fig:first}}
\subfigure[\footnotesize PID control]{\includegraphics[width=2in]{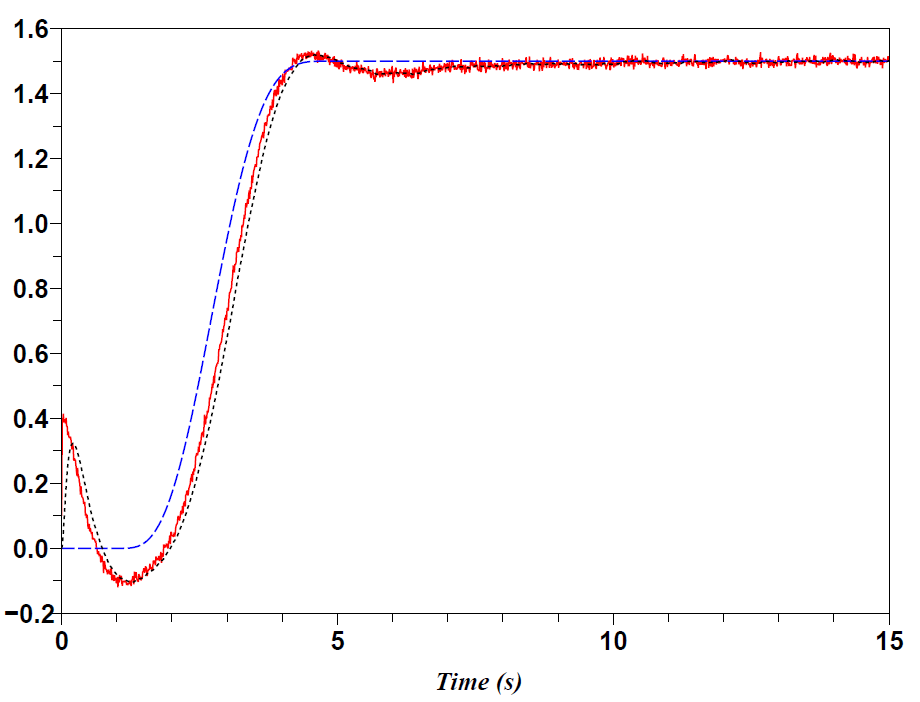}\label{fig:sec}}
\subfigure[\footnotesize i-PI control]{\includegraphics[width=2in]{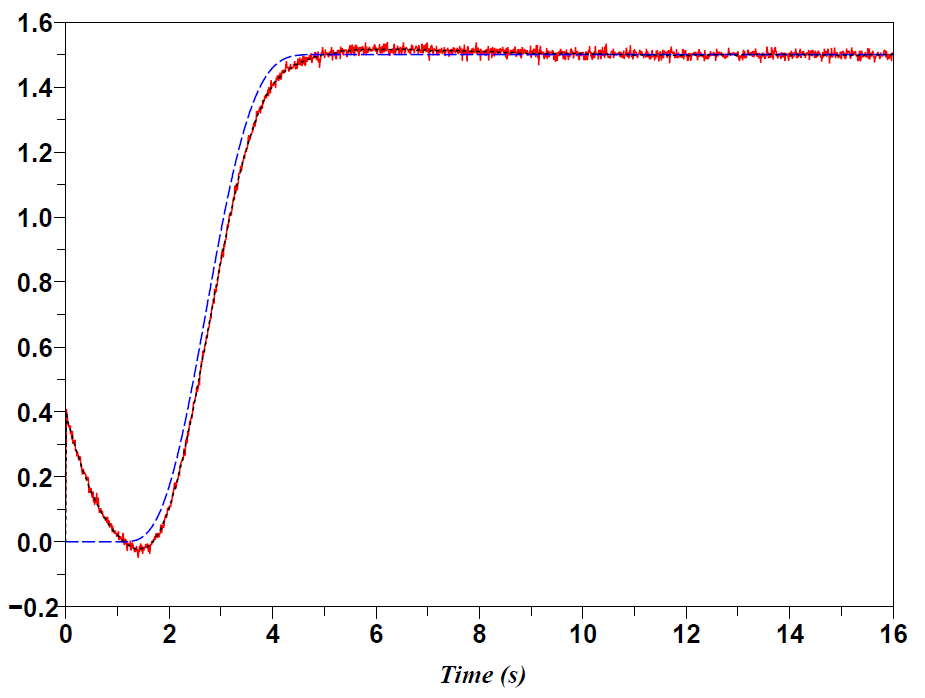}\label{fig:third}}
\subfigure[\footnotesize PID control]{\includegraphics[width=2in]{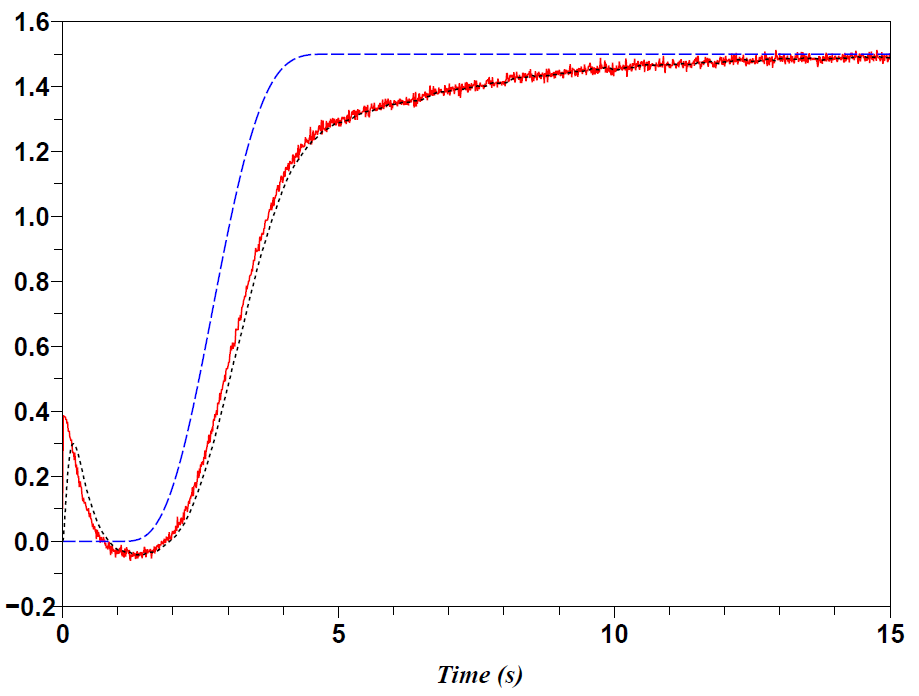}\label{fig:fourth}}
\caption{Stable linear monovariable system (Output (--); reference (- -); denoised output (. .)).}
\end{figure}

\vspace{5cm}
\section{Control of non-minimum phase systems}

We explain in this section, how to derive the model-free control law (\ref{eq:ipi}) in order to stabilize and guarantee certain
performances for non-minimum phase systems. We will show that the proposed control law is also robust to disturbances and switched models.

\subsection{Discrete model-free control law for non-minimum phase systems}

Firstly, consider the discretized model-free control law, which is typically used for a digital implementation.

\begin{mydef} \cite{Michel} \label{mydef_iPI_discret}
For any discrete moment $t_k, \, k \in \mathbb{N}$, one defines the discrete controller i-PI.
\begin{equation}\label{eq:iPI_discret_eq}
u_k =  u_{k-1} - \frac{1}{\alpha}\left(\left. {y}^{(n)}\right|_{k-1} - \left. {y}^{(n)\ast} \right|_k \right) + \left. \mathcal{C}(\varepsilon) \right|_k
\end{equation}
where
\begin{itemize}
\renewcommand{\labelitemi}{$\bullet$}
\item $y^\ast$ is the output reference trajectory;
\item $\varepsilon = y^\ast - y$ is the tracking error;
\item $\mathcal{C}$ is a usual corrector PI where $K_P$, $K_I$ are the usual tuning gains.
\end{itemize}
The discrete intelligent controller is also called discrete model-free control law or discrete model-free law.
\end{mydef}

Non-minimum phase systems are characterized by negative zero(s). Such zero can be approximated by a delay
since $e^{-Ts} \approx 1 - T s$ using a Taylor expansion. To compensate the effect of the delay, that may destabilize the control,
we take the derivative of the output $y$ instead of the output $y$ to create the measurement feedback. This way allows to anticipate
the variations of $y$ and finally cancel the disturbances associated to the presence of the delay.

We define consequently the i*-PI controller for non-minimum phase systems.

\begin{mydef}  \label{mydef_iPI_discret}
For any discrete moment $t_k, \, k \in \mathbb{N}$, one defines the discrete controller i*-PI for non-minimum phase systems.
$\lambda$ and $\delta_j$ are real coefficients.
\begin{equation}\label{eq:iPI_discret_nm_eq}
u_k =  \mathcal{G}(\varepsilon) \left\{ u_{k-1} - \sum_{j=1}^n \delta_j \left( \lambda \left. {y}^{(j)}\right|_{k-1} - \left. {y}^{(j)\ast} \right|_k \right) \right\}
\end{equation}
where
\begin{itemize}
\renewcommand{\labelitemi}{$\bullet$}
\item $y^\ast$ is the output reference trajectory;
\item $\varepsilon = y^\ast - y$ is the tracking error;
\item $\mathcal{G(\varepsilon)}$ is called a gain function and is either a pure gain or an
integrator\footnote{Depending on the application, a pure gain can be enough to ensure good tracking performances.}.
\end{itemize}
The discrete intelligent controller is also called discrete NM-model-free control law or discrete NM-model-free law.
\end{mydef}

\noindent
Practically, simulations show that $n = 2$ is sufficient\footnote{The possibility of reducing $n$ will be studied in a future work.} to ensure at least the stability of the model-free control closed-loop. Therefore, (\ref{eq:iPI_discret_nm_eq}) is written :

\begin{equation}\label{eq:iPI_discret_nm_eq_2}
u_k =  \mathcal{G}(\varepsilon) \left\{ u_{k-1} - \delta_2  \left( \lambda \left. \frac{\diff^2 y}{\diff t^2}  \right|_{k-1} - \left. \frac{\diff^2 y^{\ast}}{\diff t^2}  \right|_{k}  \right) - \delta_1 \left( \lambda \left. \frac{\diff y}{\diff t}  \right|_{k-1} - \left. \frac{\diff y^{\ast}}{\diff t} \right|_{k}  \right)  \right\}
\end{equation}

\noindent
For the following applications, we choose the gain function as an integrator, with a $K_i$ constant, such that :

\begin{equation}\label{eq:iPI_discret_nm_eq_2}
\mathcal{G}(\varepsilon) = K_i \int_0^t \varepsilon \, \diff t
\end{equation}

\subsection{Applications}

Consider the systems $\Sigma_1$, $\Sigma_2$, $\Sigma_3$ and $\Sigma_4$, which are minimum and non-minimum phase systems, and which are described respectively by the state-space representations :

\begin{equation}
\Sigma_1 := \left\{ \begin{array}{l}
\dot{x} = \begin{pmatrix}
            0 & -1000 \\
         100000 & -5000
          \end{pmatrix} x  +
          \begin{pmatrix}
            2.10^4 \\
             0
          \end{pmatrix} u \\
 y = \begin{pmatrix}
        -10 & 1
      \end{pmatrix} x
      \end{array} \right.
\end{equation}

\vspace{0.5cm}
\begin{equation}
\Sigma_2 := \left\{ \begin{array}{l}
\dot{x} = \begin{pmatrix}
            0 & -900 \\
            80000 & -3500
          \end{pmatrix} x  +
          \begin{pmatrix}
            2.10^4 \\
             0
          \end{pmatrix} u \\
 y = \begin{pmatrix}
        -13 & 1
      \end{pmatrix} x
      \end{array} \right.
\end{equation}

\vspace{0.5cm}
\begin{equation}
\Sigma_3 := \left\{ \begin{array}{l}
\dot{x} = \begin{pmatrix}
            0 & -900 \\
            80000 & -3500
          \end{pmatrix} x  +
          \begin{pmatrix}
            2.10^4 \\
             0
          \end{pmatrix} u \\
 y = \begin{pmatrix}
        +13 & 1
      \end{pmatrix} x
      \end{array} \right.
\end{equation}

\vspace{0.5cm}
\begin{equation}
\Sigma_4 := \left\{ \begin{array}{l}
\dot{x} = \begin{pmatrix}
            0 & -400  \\
            70000 & -1500
          \end{pmatrix} x  +
          \begin{pmatrix}
            2.10^4 \\
             0
          \end{pmatrix} u \\
 y = \begin{pmatrix}
        +5 & 1
      \end{pmatrix} x
      \end{array} \right.
\end{equation}

\begin{figure}[!h]
\centering
\includegraphics[width=12cm]{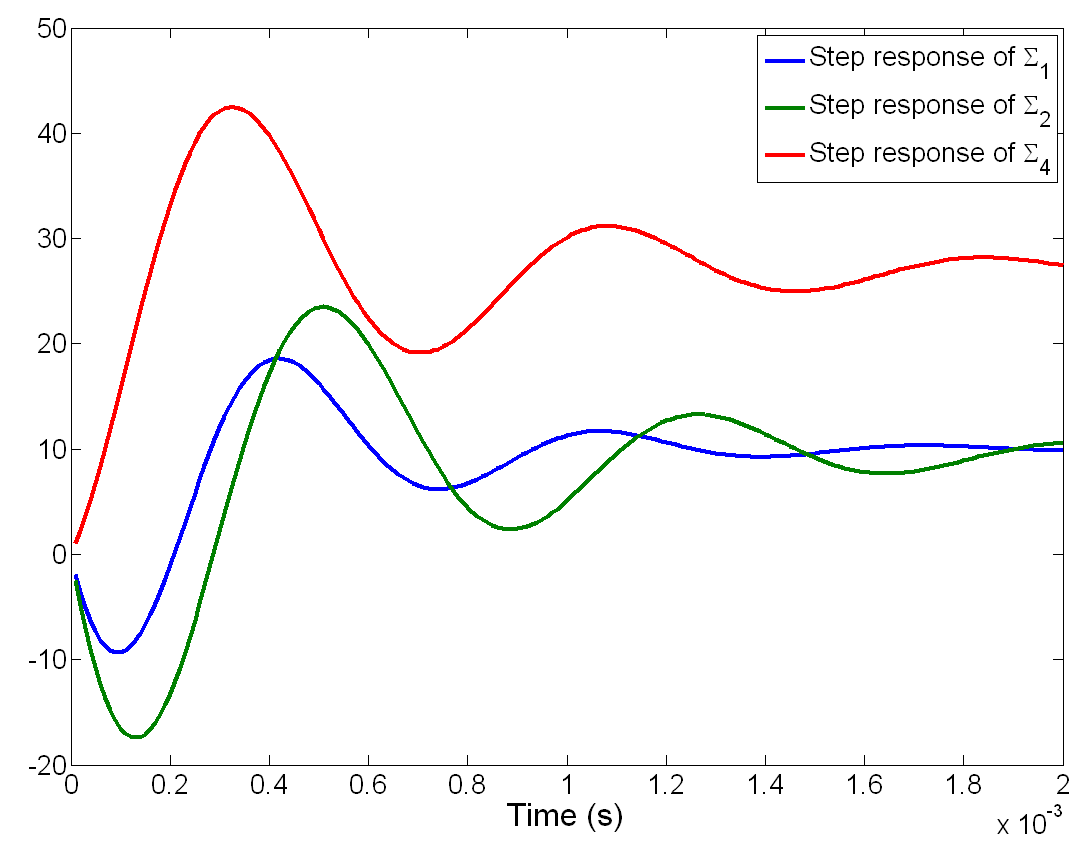}
\caption{Step responses of the systems $\Sigma_1, \, \Sigma_2$ and $\Sigma_4$.}
\label{fig:NM_Fig_1}
\end{figure}

The unitary step response of these systems is presented Fig. \ref{fig:NM_Fig_1}.

\begin{figure}[!h]
\centering
\subfigure[\footnotesize Transient with stabilization.]{\includegraphics[width=11cm]{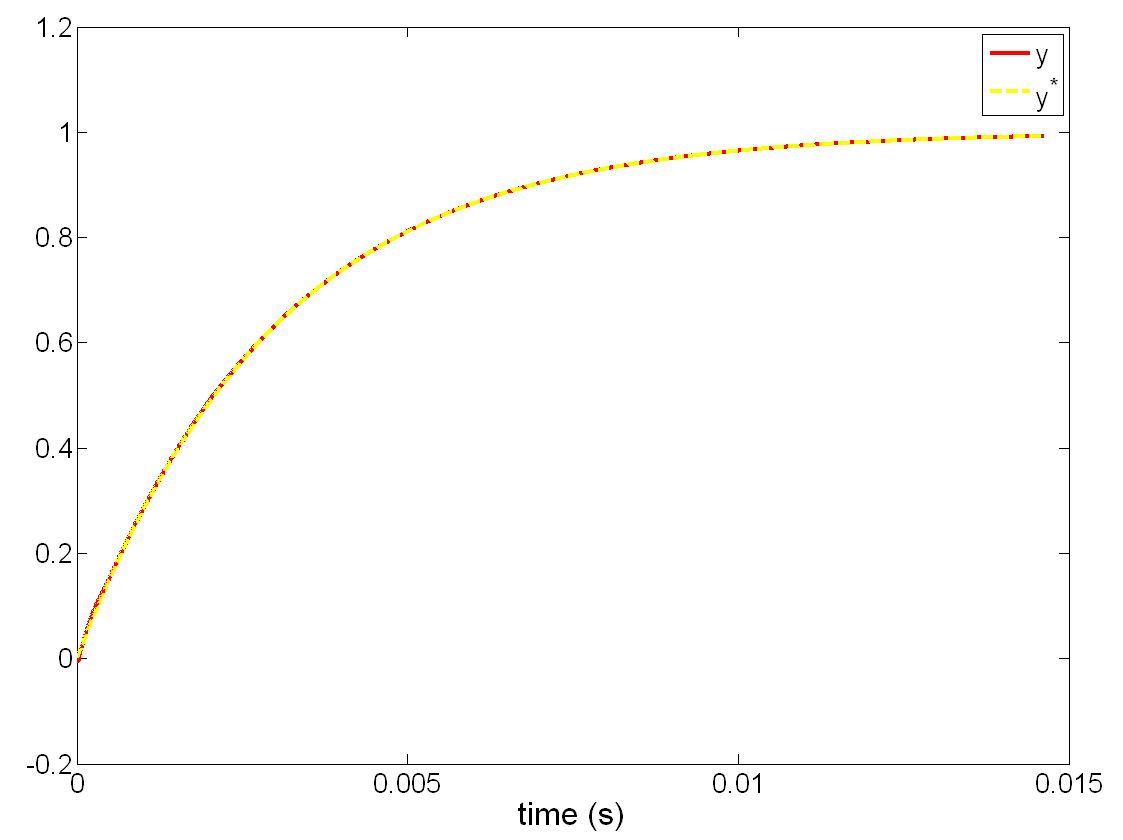}\label{fig:NM_Fig_2}}
\subfigure[\footnotesize Focus on the beginning of the transient.]{\includegraphics[width=11cm]{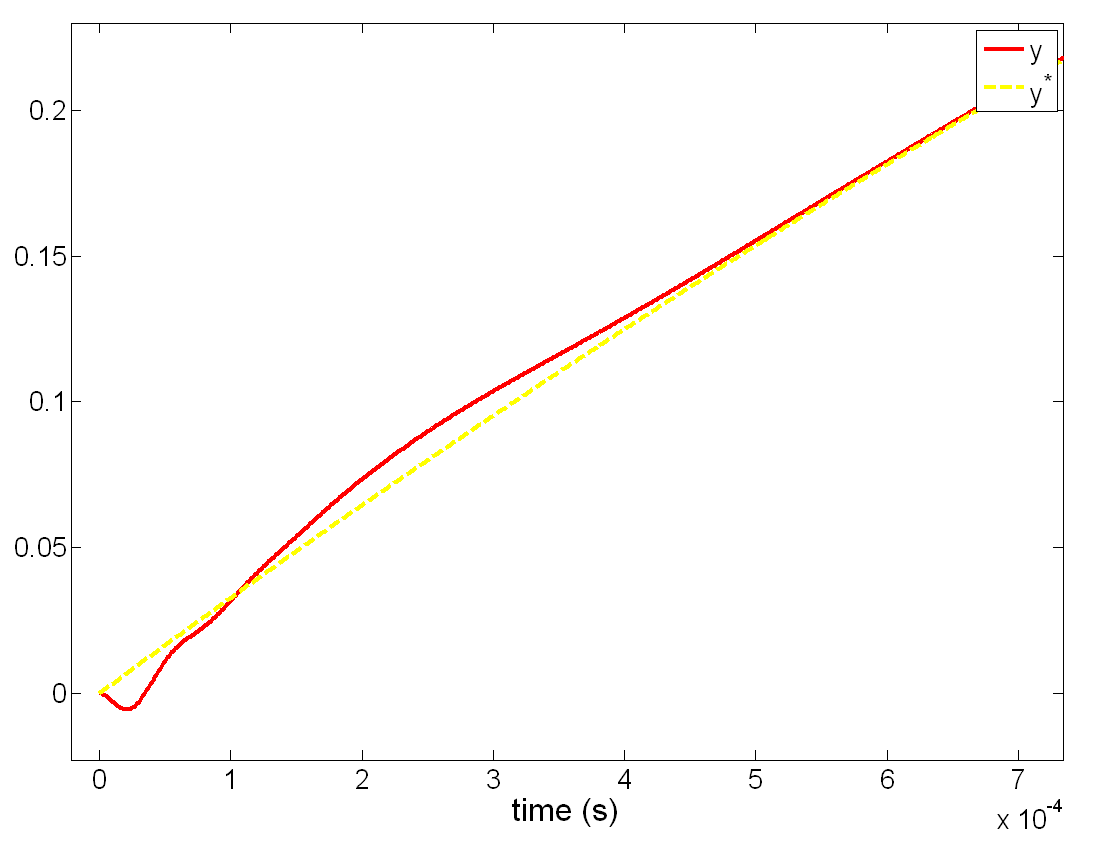}\label{fig:NM_Fig_3}}
\caption{Tracking of an exponential reference. }
\label{fig:NM_Fig_2_}
\end{figure}

\begin{figure}[!h]
\centering
\includegraphics[width=11cm]{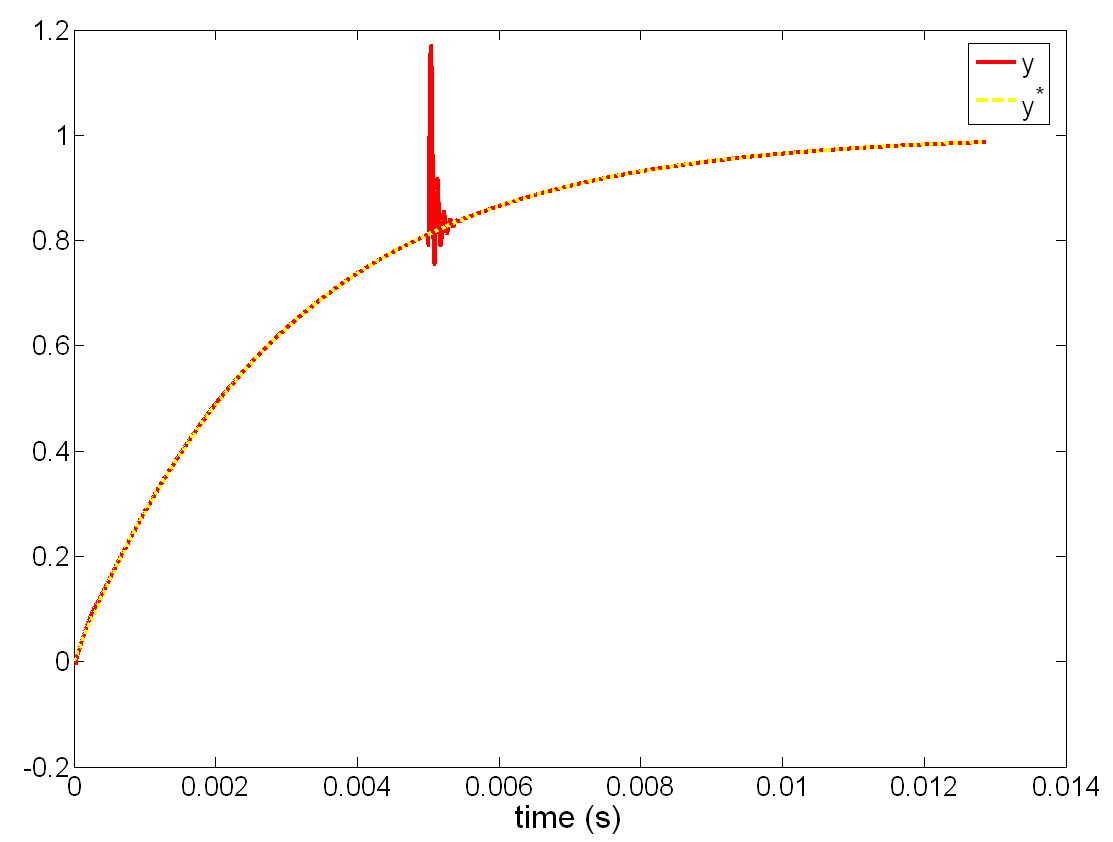}
\caption{Tracking of an exponential reference; $\Sigma_1$ switches to $\Sigma_2$ at $t = 5$ ms.}
\label{fig:NM_Fig_4}
\vspace{0.5cm}
\includegraphics[width=11cm]{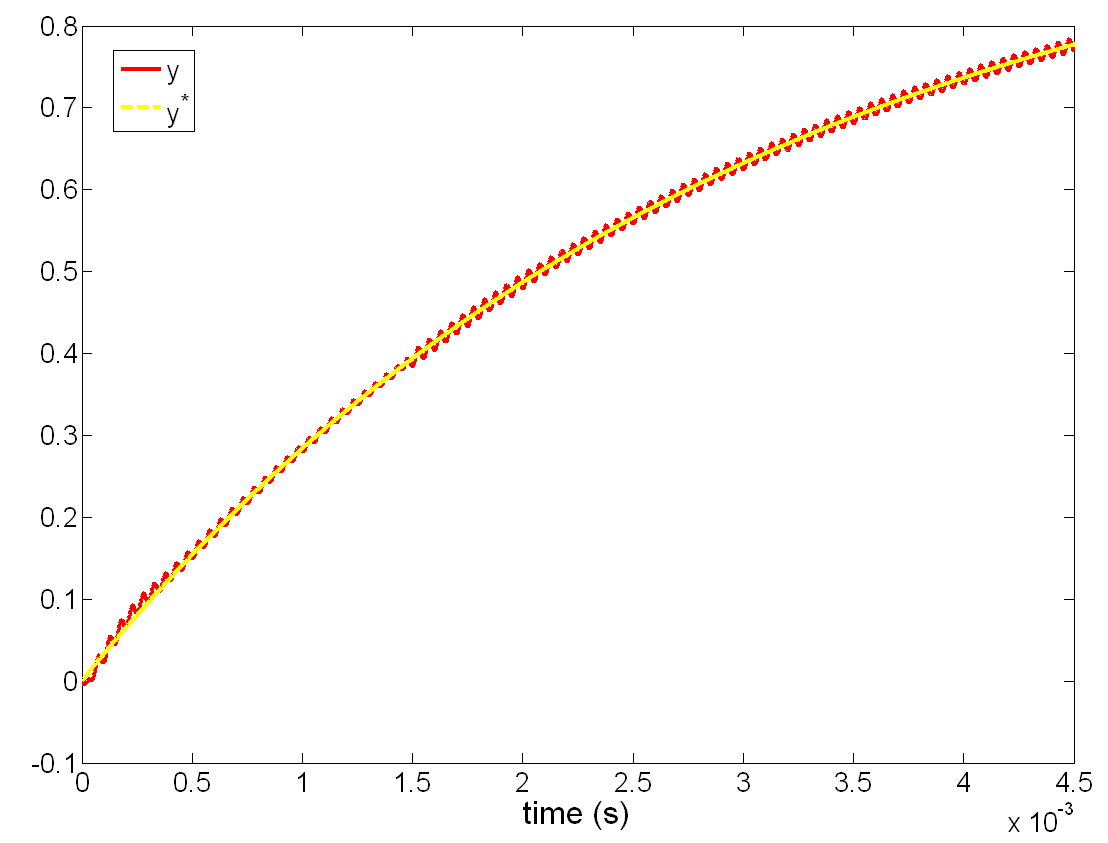}
\caption{Tracking of an exponential reference with a sinusoidal disturbance added on $u$ such that $\tilde{u} = 5 \cos(\frac{2 \pi}{5.10^{-5}}t)$.}
\label{fig:NM_Fig_5}
\end{figure}

\begin{figure}[!h]
\centering
\includegraphics[width=11cm]{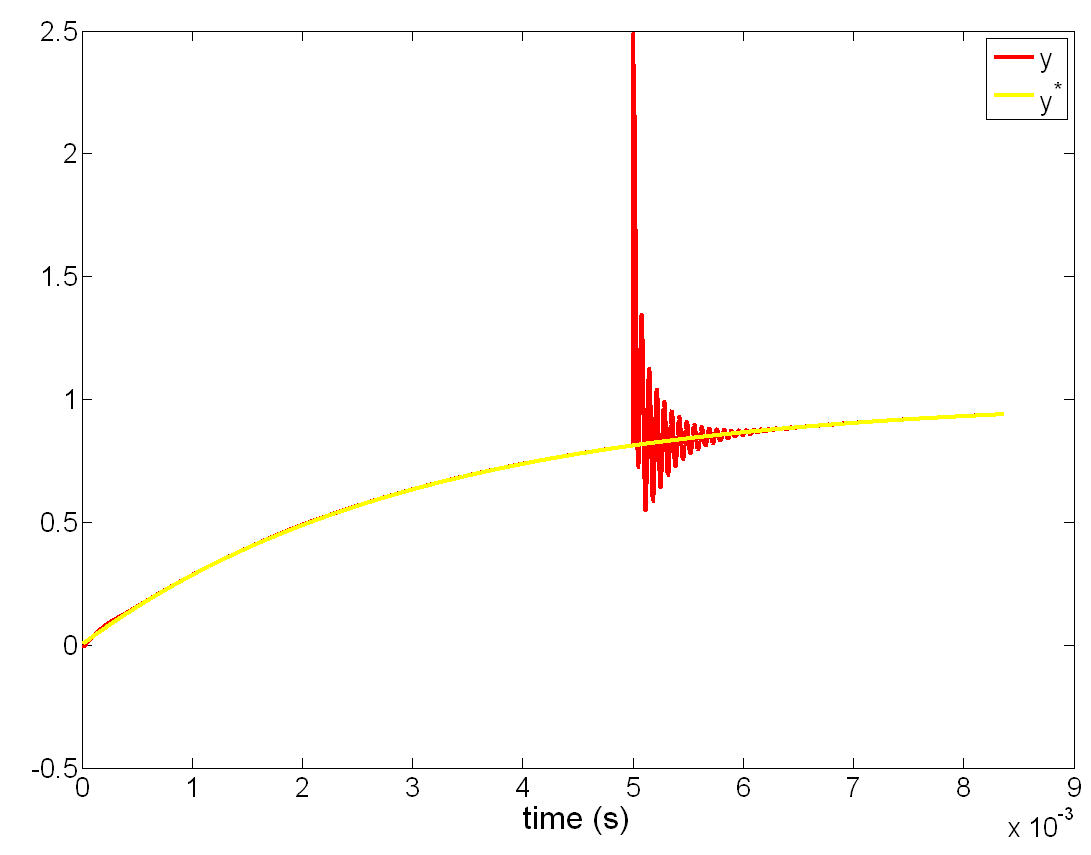}
\caption{Tracking of an exponential reference; $\Sigma_1$ switches to $\Sigma_3$ at $t = 5$ ms.}
\label{fig:NM_Fig_6}
\vspace{0.5cm}
\includegraphics[width=11cm]{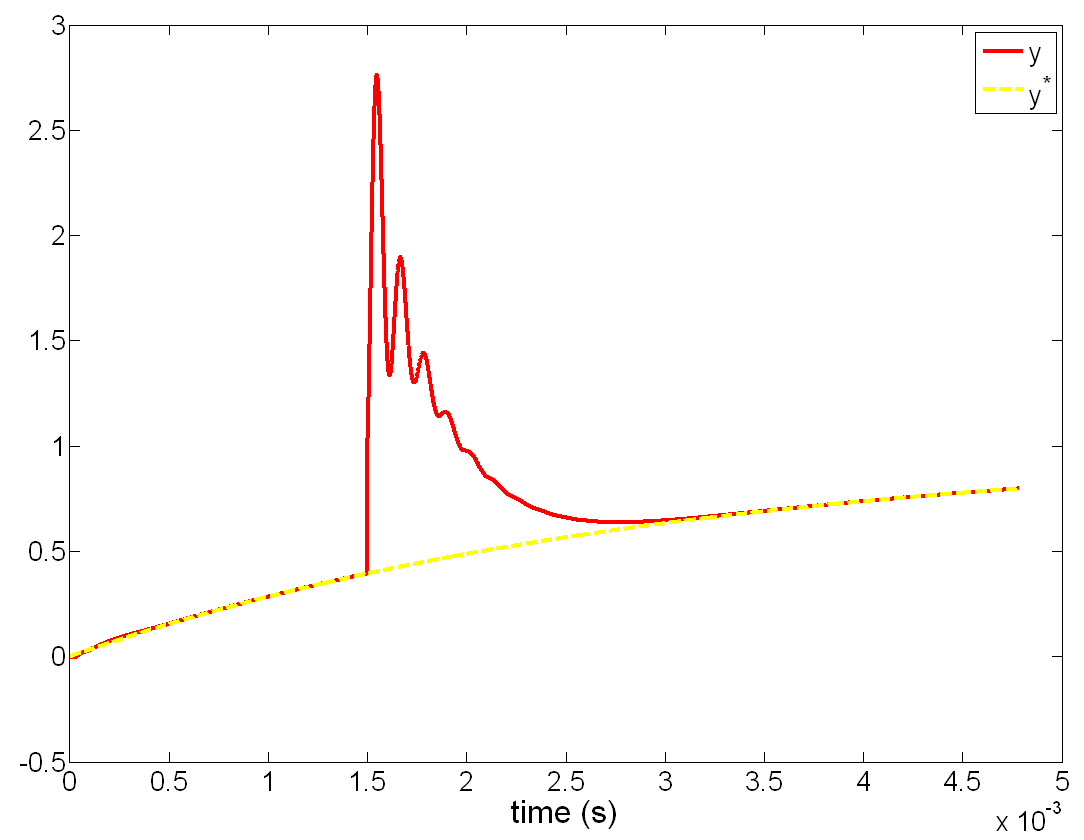}
\caption{Tracking of an exponential reference; $\Sigma_1$ switches to $\Sigma_4$ at $t = 1.5$ ms.}
\label{fig:NM_Fig_7}
\end{figure}

\begin{figure}[!h]
\centering
\includegraphics[width=11cm]{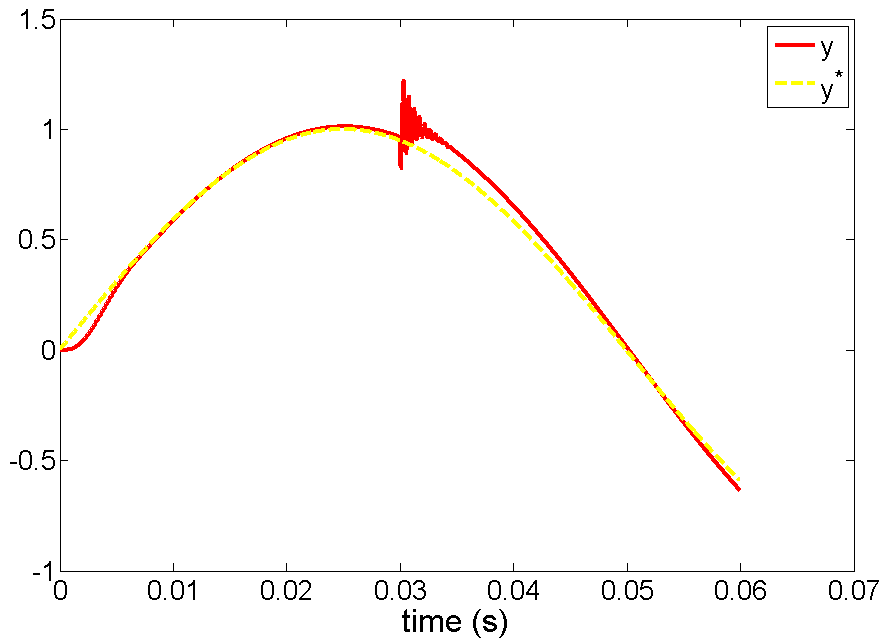}
\caption{Tracking of a sinusoidal reference; $\Sigma_1$ switches to $\Sigma_2$ at $t = 3$ ms.}
\label{fig:NM_Fig_8}
\vspace{0.5cm}
\includegraphics[width=11cm]{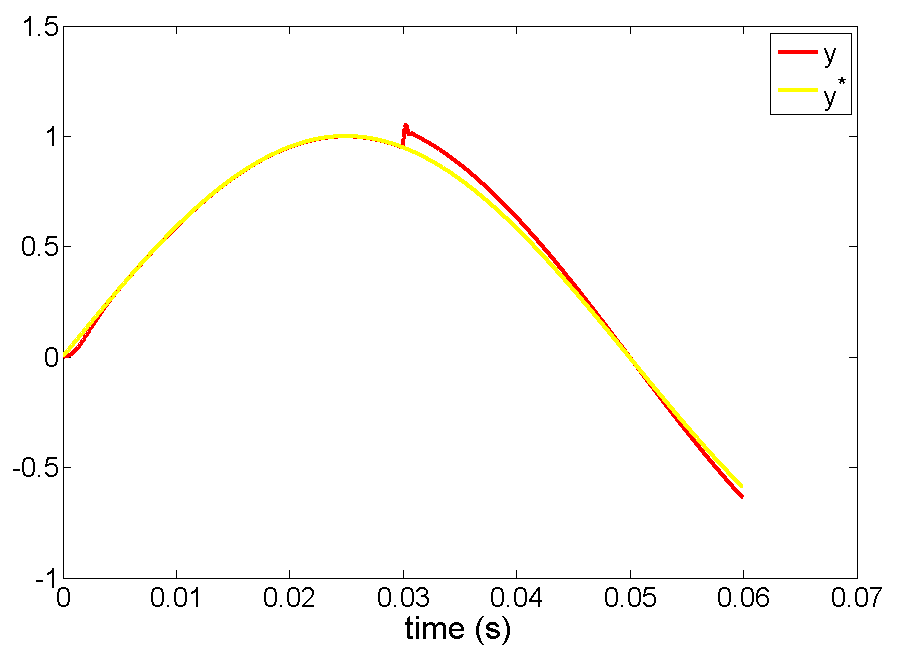}
\caption{Tracking of a sinusoidal reference; $\Sigma_1$ switches to $\Sigma_3$ at $t = 3$ ms.}
\label{fig:NM_Fig_9}
\end{figure}

The following figures present some preliminary examples of the application of the i*-PI control.
Figure \ref{fig:NM_Fig_2_} presents the tracking of an exponential reference for the system $\Sigma_1$ (with a focus on the beginning of the transient). Figures \ref{fig:NM_Fig_4} and \ref{fig:NM_Fig_5} show the response $y$ of the controlled system  when respectively $\Sigma_1$ switches to $\Sigma_2$ and with the addition of a sinusoidal disturbance on the variable $u$. Figures \ref{fig:NM_Fig_6} and \ref{fig:NM_Fig_7} present the control of switched systems; in particular the commutation from a non-minimum phase system to a minimum phase system. Figures \ref{fig:NM_Fig_8} and \ref{fig:NM_Fig_9} present the tracking of a sinusoidal reference when systems switch. The case where $\Sigma_1$ switches to $\Sigma_3$ has been already studied in \cite{Michel_microgrids}. We investigated the application of the model-free control in a microgrid environment under load / transfer function changes. These changes imply substantial modifications of the controlled models.  

%
%
%
\newpage
\section{Concluding remarks}

We presented some simulation results that confirm the fact that the NM-model-free control or i*-PI controller, designed for the control of non-minimum phase systems, ensure good tracking performances. We evaluated the performances in presence of disturbances and in the case of switched systems. In particular, the NM-model-free control is able {\it a priori} to control both minimum and non-minimum-phase systems. The proposed control law seems to have the same properties than the original model-free control \cite{esta} for which its performances have been successfully proved in simulation  when controlling switched systems (e.g. \cite{Michel_microgrids}). Further work will concern the study of the stability of the NM-model-free control method and its applications to networked systems.

\end{document}